\documentclass[12pt]{amsart}
\usepackage{amssymb}
\newtheorem{theorem}[subsection]{Theorem}
\newtheorem{proposition}[subsection]{Proposition}
\newtheorem{corollary}[subsection]{Corollary}

\theoremstyle{remark}

\numberwithin{equation}{section}

\newcommand{\End}{{\mathrm{End}}}

\def\Bbb{\mathbb}
\let\i=\iota

\def\endrk{\hbox{$|\!\!|\!\!|\!\!|\!\date!|\!\!|\!\!|$}}

\begin{document}

\title{Invariant Prolongation and Detour Complexes}
\dedicatory{Dedicated to the memory of
Thomas P. Branson (1953 - 2006)}

\author{A.~Rod Gover}
\address{Department of Mathematics, University of Auckland,
Private Bag 92019,\newline\indent Auckland, New Zealand}
\email{gover@math.auckland.ac.nz}

\newcommand{\nn}[1]{(\ref{#1})}

\def\sideremark#1{\ifvmode\leavevmode\fi\vadjust{\vbox to0pt{\vss
 \hbox to 0pt{\hskip\hsize\hskip1em
 \vbox{\hsize3cm\tiny\raggedright\pretolerance10000
 \noindent #1\hfill}\hss}\vbox to8pt{\vfil}\vss}}}%

                                                   %
\newcommand{\edz}[1]{\sideremark{#1}}

\thanks{The author would like to thank the Institute for Mathematics
  \& Its Applications, University of Minnesota, and also the
  organisers of the 2006 Summer Programme there ``Symmetries and
  overdetermined systems of partial differential equations''. Special
  thanks are due to Tom Branson, a friend, a colleague and a mentor.
  It is my understanding that his enthusiasm and work contributed
  significantly to the occurrence and organisation of the meeting. It
  is part of the tragedy of his departure that he did not get to  enjoy the
  realisation of this event.}

\begin{abstract}
  In these expository notes we draw together and develop the ideas
  behind some recent progress in two directions: the treatment of
  finite type partial differential operators by prolongation, and a
  class of differential complexes known as detour complexes. This
  elaborates on a lecture given at the IMA Summer Programme 
  ``Symmetries and overdetermined systems of partial differential
  equations''.
\end{abstract}

\maketitle

\section{Introduction}

Differential complexes capture integrability conditions for linear
partial differential operators. When they involve sequences of natural
operators on geometric structures, such complexes typically encode deep
geometric information about the structure. In these expository notes
we survey some recent results concerning the construction of such
complexes. A particular focus is the construction via prolonged
differential systems and the links to the treatment of finite type
differential operators. The broad issues arise in many contexts but we
take our motivation here from conformal and Riemannian geometry.

Elliptic differential operators with good conformal behaviour play a
special role in geometric analysis on Riemannian manifolds.  The
``Yamabe Problem'' of finding, via conformal rescaling, constant
scalar curvature metrics on compact manifolds is a case in point.
This exploits heavily the conformal Laplacian, since it controls the
conformal variation of the scalar curvature (see \cite{L&P} for an
overview and references). The higher order conformal Laplacian
operators of Paneitz, Riegert, Graham et al.\ \cite{GJMS} (the {\em
  GJMS operators}) have been brought to bear on related problems by
Branson, Chang, Yang and others \cite{Brsharp,aliceC-EMS,GrZ}.
However for many important tensor or spinor fields there is no natural
conformally invariant elliptic operator (taking values in an
irreducible bundle) available.  This is true even in the conformally
flat case, and indeed this claim follows easily from the
classification of conformal differential operators on the sphere
\cite{BoeColl}.  From this it is clear that for many bundles on the
sphere the analogue, or replacement, for an elliptic operator is an
elliptic complex of conformally invariant differential operators.
Ignoring the issue of ellipticity, the requirement that a
sequence of operators be both conformally invariant and form a complex
is already severe, and so the construction of such complexes is a delicate
matter.

On the conformal sphere there is a class of conformal complexes known
as Bernstein-Gelfand-Gelfand (BGG) complexes.  Equivalents of these
were first constructed in the representation theory of (generalised)
Verma modules (see \cite{Lepow} and references therein).  From there
we see, more generally, versions of these complexes for a large class
of homogeneous structures in the category of so-called parabolic
geometries, that is manifolds of the form $G/P$ where $G$ is a
semisimple Lie group and $P$ a parabolic subgroup. Although these
early constructions were motivated by geometric questions, the entry
of the unifying picture of these complexes into the mainstream of
Differential Geometry was pioneered by Eastwood, Baston and
collaborators \cite{ER,Baston}. In particular they developed
techniques for constructing sequences of differential operators on
general curved backgrounds that generalise the BGG sequences. Through
work over the years the nature of many of the operators involved has
become quite well developed \cite{Goquart,CDS}, and there has been
spectacular recent progress in general approaches to constructing
these complexes \cite{CSS,CD}. However, as we shall discuss in Section
\ref{detseqs}, these sequences are in general not complexes on curved
manifolds.

In the BGG sequences each operator may be viewed as an integrability
condition for the preceding operator. The operators encode in part
the geometry and it seems likely that in general, at least for
conformal structures, one obtains a complex only in the conformally
flat setting.  In work with Tom Branson the ambient metric of
Fefferman-Graham \cite{FeffGrast} and then variational ideas were used
to give two constructions \cite{BrGodeRham,BrGodefdet} of what we have
termed ``detour complexes''. These are related to BGG complexes, but
involve fewer operators, weaker integrability conditions and yield
complexes in curved settings. The central results and ideas (with some
extensions) are reviewed and developed in Sections \ref{detseqs},
\ref{curveddets} and \ref{var}. The complexes obtained there have
interesting geometric applications and interpretations, but the
constructions involved do not suggest how a full theory of such
complexes could be developed. We should note that considerably
earlier, following related constructions of Calabi and
Gasqui-Goldschmidt, Gasqui constructed \cite{Gasqui} a three step
elliptic sequence which forms a complex if and only if the manifold is
Einstein. This was studied further in \cite{BGilkeyZ}. In our current
language this is a detour complex, exactly along lines we discuss, and
may be viewed as a Riemannian analogue of the conformal complex of
\cite{BrGodefdet} as in Theorem \ref{obflat}.

Beginning at Section \ref{YM} we start the development of a model for
a rather general construction principle for a class of detour
complexes on pseudo-Riemannian manifolds. This is based on recent work
with Petr Somberg and Vladimir Sou\v cek.  We treat two examples, one
in some detail, and the complexes we obtain are elliptic in the case
that the background structure has Riemannian signature. In dimension 4
the complexes are conformally invariant. There are two main tools
involved in the constructions. The first is general sequence which
yields a detour complex for every Yang-Mills connection. The second
ingredient is the use of prolonged differential systems.

The use of prolonged systems is related to the following problem:
Given a suitable (i.e. finite type, as explained below) linear partial
differential operator $D$, can one construct an equivalent first order
prolonged system that is actually a vector bundle connection?  Here
there are close connections with the recent work \cite{BCEG} with
Branson, \v Cap and Eastwood (for this, and further discussion of BGG 
complexes, see also the proceedings of A.\ \v{C}ap \cite{CIMA}). In fact we want more.  The connection
should be invariant in the sense that it should share the symmetry and
invariance properties of $D$.

\section{Invariant tractor connections for finite type PDE}\label{prolongs}

First we introduce some notation to simplify later presentations.  We
shall write $\Lambda^1$ for the cotangent bundle, and $\Lambda^p$ for
its $p^{\rm th}$ exterior power. Thus the trivial bundle (with fibre
$\mathbb{R}$) may be denoted $\Lambda^0$ although, for reasons that
will later become obvious, we will also write $\mathcal{E}$ for this.
For simplicity the $k^{\rm th}$ symmetric powers of $\Lambda^1$ will
be denoted $S^k$ (with the metric trace-free part of this denoted $S^k_0$);
albeit that this introduces the redundancy that $S^1=\Lambda^1$. In
all cases we will use the same notation for a bundle as for its
section spaces. All structures will be assumed smooth and we restrict
to differential operators that take smooth sections to smooth
sections.

 The prolongations of a $k^{\rm th}$-order semilinear differential
 operator $D:E\to F$, between vector bundles, are constructed from its
 leading symbol $\sigma(D):S^k  \otimes E\to F$.
At a point of $M$, denoting by $K^0$ the kernel of
 $\sigma(D)$, the spaces $V_i=(S^i\otimes E)\cap
 (S^{i-k}\otimes K^0)$, $i\geq k$, capture spaces of new
 variables to be introduced and the system closes up and is said to be
 of {\em finite type} if $V_i=0$ for sufficiently large $i$. 
 For
 example, on a Riemmanian manifold $(M^n,g)$ a vector field $k$ is an
 infinitesimal isometry, a so-called Killing vector field, if Lie
 differentiation along its flow preserves the metric $g$, that is
 $\mathcal{L}_k g=0$.  Rewriting this in terms of the Levi-Civita
 connection $\nabla$ (i.e. the unique torsion free connection on $TM$
 preserving the metric) the Killing equation is seen to be a system of
 $n(n+1)/2$ equations in $n$ unkowns, viz.\ $\nabla_a k_b+\nabla_b
 k_a=0$, where we have used an obvious abstract index notation and
 also we have used the metric to identify $k$ with the 1-form
 $g(k,~)$. In this example 
 $K^0$ is evidently $\Lambda^2$, the space of 2-forms, and we may rewrite
 the system as $\nabla_a k_b=\mu_{ab}$ where $\mu\in \Lambda^2$ (that
 is $\mu_{ab}=-\mu_{ba}$). Since $(S^2\otimes \Lambda^1)\cap
 (\Lambda^1\otimes \Lambda^2)=0$, when we differentiate and compute
 consequences the system closes up algebraically after just one step.
 We obtain a prolonged system which is the equation of parallel
 transport for a connection:

\vspace{-2mm}

\begin{equation}\label{projex}
\nabla^D_a\left(\begin{array}{c}k_b\\\mu_{bc}\end{array}\right):=
\left(\begin{array}{c}\nabla_a k_b-\mu_{ab}\\\nabla_a \mu_{bc}-R_{bc}{}^d{}_a k_d\end{array}\right)=0,
\end{equation}

\vspace{-.5mm}

\noindent where $R$ is the curvature of $\nabla$.  We will regard this
connection as being {\em equivalent} to the original equation, since solutions
of the original equation are in 1-1 correspondence with sections of
$\mathbb{T}:=\Lambda^1\oplus \Lambda^2 $ that are parallel for
$\nabla^D$.  It follows that the original equation has a solution
space of dimension at most rank$(\mathbb{T})=n(n+1)/2$.  The curvature
of $\nabla^D$ evidently obstructs solutions and, in particular, the
maximal number of solutions is achieved only if the connection
$\nabla^D$ is flat.

\vspace{1.3mm}

Although this example is very simple it already brings us to several
of the essential issues. The original equation $\mathcal{L}_k g=0$
arose in a Riemannian setting.  By the
explicit formula given it is clear that that the connection $\nabla^D$ is
naturally invariant for Riemannian structures.  However comparing with
section 3 in \cite{BEGo}, or \cite{CapGoluminy}, one identifies this
connection as an obvious curvature modification of the normal tractor
connection on projective structures. The Levi-Civita connection
$\nabla$ in \nn{projex} may be viewed as an affine connection, so this
makes sense. If $\nabla$ is any torsion-free connection on $\Lambda^1$
then note that the equation $\nabla_a k_b+\nabla_b k_a=0$ is invariant
under the transformations
$$
k_b\mapsto \widehat{k}_b = e^{2\omega}k_b \quad \mbox{and} \quad \nabla\mapsto \widehat{\nabla}
$$ 
where, as an operator on 1-forms, 
$$ \widehat{\nabla}_a u_b= \nabla_a u_b -\Upsilon_a u_b -\Upsilon_b
u_a \quad \mbox{with} \quad \Upsilon =d\omega~.
$$ These show that the equation is in fact {\em projectively
  invariant}, it is well defined on manifolds having only an
equivalence class of torsion-free connections, where the equivalence
relation is given by the transformations indicated, i.e. on projective
manifolds. This leads us to question whether the tractor connection
$\nabla^D$ in \nn{projex} shares this property. It does.  It is readily
verified directly that it is also projectively invariant; it differs
from the invariant normal projective connection in \cite{BEGo,CapGoluminy}
by an action of the
(projectively invariant) curvature of the normal connection. In this
case it is straightforward to see that this was inevitable.

In general prolonged systems are more complicated.  In \cite{BCEG}
Kostant's algebraic Hodge theory \cite{Kost61} led to an explicit and
uniform treatment of prolongations for a large class of overdetermined
PDE. One of the simplest classes of examples is the set of equations
controlling conformal Killing forms and via explicit calculations in
terms of the Levi-Civita connection, prolonged systems for these were
earlier calculated by \cite{semm}. However in neither of these works
was the issue of invariance discussed. The equation for a conformal
Killing $p$-form $\kappa$ can be given as follows: for any tangent
vector field $u$ we have
$$
\nabla_u \kappa = \varepsilon(u)\tau + \iota(u) \rho
$$ where, on the right-hand side $\tau$ is a $(p-1)$-form, $\rho$ is a
$(p+1)$-form, and $\varepsilon(u)$ and $\iota(u)$ indicate,
respectively, the exterior multiplication and (its formal adjoint) the
interior multiplication of $g(u,~)$. An important property of the
conformal Killing equation is that it is conformally invariant (where
we require the $p$-form $\kappa$ to have conformal weight $p+1$).
This can be phrased in similar terms to the projective invariance of
the equation described above, but alternatively it simply means the
equation descends to a well defined equation on manifolds where we do
not have a metric, but rather only an equivalence class of conformally
related metrics (this is a conformal structure).  So it is natural to ask
whether there is an equivalent prolonged system that may be realised
as a conformally invariant connection. Using the framework of tractor calculus,
this is answered in the affirmative in \cite{GoSil}, where also the
invariant connection is related to the normal conformal tractor
connection.

The leading symbol determines whether or not an equation is of finite
type; an operator $D:E\to F$ is of finite type if and only if its
(complex) characteristic variety is empty \cite{spencer}.  One may ask
whether any finite type overdetermined linear differential operator is
equivalent, in the sense of prolongations, to a connection on vector
bundle where the connection is ``as invariant'' as the original
operator, that is it shares the same symmetries. This is obviously an
important question in its right. For example if for some PDE one finds
a connection on prolonged system which depends only on operator $D$
and no other choices, then it is invariant in this sense and its
curvature is really a geometric invariant of the equation $D$. In
general this is probably more than one can hope for. It is probably
only reasonable to hope to find a canonical connection after some well
defined choices controlled by representation theory on finite
dimensional vector spaces that are associated with the original equation.  It
turns out such questions are also important for the construction of
certain natural differential complexes.

\section{Detour complexes}\label{detoursec}

\newcommand{\si}{\sigma}
\newcommand{\ce}{\mathcal{E}}
\newcommand{\bg}{\mbox{\boldmath{$ g$}}}
\newcommand{\bT}{\mathcal{T}}
\newcommand{\bD}{\Bbb D}    
\newcommand{\h}{\mbox{\boldmath{$ h$}}}
\newcommand{\Om}{\Omega}
\newcommand{\cq}{\mathcal{Q}}
\newcommand{\Ric}{\operatorname{Ric}}
\newcommand{\Sc}{\operatorname{Sc}}
\newcommand{\id}{\operatorname{id}}

We we will specialise our discussions here to the setting of oriented
(pseudo-)Riemannian structures and conformal geometries. The restriction to
oriented structures is just to simplify statements and is not
otherwise required.  It should also be pointed out that similar
ideas apply in many other settings including, for example, CR
geometry.

Recall that a {\em conformal manifold} of signature $(p,q)$ on $M$ is
a smooth ray subbundle $\cq\subset S^2T^*M$ whose fibre over $x$
consists of conformally related signature-$(p,q)$ metrics at the point
$x$. Sections of $\cq$ are metrics $g$ on $M$. So we may equivalently
view the conformal structure as the equivalence class $[g]$ of these
conformally related metrics.  The principal bundle $\pi:\cq\to M$ has
structure group $\Bbb R_+$, and so each representation ${\Bbb R}_+ \ni
x\mapsto x^{-w/2}\in {\rm End}(\Bbb R)$ induces a natural line bundle
on $ (M,[g])$ that we term the conformal density bundle and denote
$\ce[w]$.  As usual we use the same notation for its section space.

We write $\bg$ for the {\em conformal metric}, that is the
tautological section of $S^2[2]:=S^2 \otimes \ce[2]$ determined by the
conformal structure. This will be used to identify $T$ with
$\Lambda^1[2]$.  For many calculations we will use abstract indices in
an obvious way.  Given a choice of metric $ g$ from the conformal
class, we write $ \nabla$ for the corresponding Levi-Civita
connection. With these conventions the Laplacian $ \Delta$ is given by
$\Delta=\bg^{ab}\nabla_a\nabla_b= \nabla^b\nabla_b\,$.  Note $\ce[w]$
is trivialised by a choice of metric $g$ from the conformal class, and
we write $\nabla$ for the connection corresponding to this
trivialisation.  It follows immediately that (the coupled) $ \nabla_a$
preserves the conformal metric.

Since the Levi-Civita connection is torsion-free, its curvature
$R_{ab}{}^c{}_d$ (the Riemannian curvature) is given by $
[\nabla_a,\nabla_b]v^c=R_{ab}{}^c{}_dv^d $ ($[\cdot,\cdot]$ indicates
the commutator bracket).  This can be decomposed into the totally
trace-free Weyl curvature $C_{abcd}$ and a remaining part described by
the symmetric {\em Schouten tensor} $P_{ab}$, according to $
R_{abcd}=C_{abcd}+2\bg_{c[a}P_{b]d}+2\bg_{d[b}P_{a]c}, $ where
$[\cdots]$ indicates antisymmetrisation over the enclosed indices.
The Schouten tensor is a trace modification of the Ricci tensor
$\Ric_{ab}$ and vice versa: $\Ric_{ab}=(n-2)P_{ab}+J\bg_{ab}$, where
we write $ J$ for the trace $ P_a{}^{a}$ of $ P$. The Weyl curvature
is conformally invariant.

\newcommand{\cB}{\mathcal{B}}
\subsection{Detour sequences}  \label{detseqs}
In conformal geometry the de Rham complex is a prototype for a 
class of sequences of bundles and conformally invariant differential
operators, each of the form 
$$
\cB^0\to \cB^1\to \cdots \to \cB^n 
$$ where the vector bundles $ \cB^i$ are irreducible tensor-spinor
bundles. For example, on the sphere with its standard conformal
structure we have the following: these are complexes and there is one
such complex for each irreducible module $\mathbb{V}$ for the group
$G=SO(n+1,1)$ of conformal motions; the space of solutions of the
first (finite type) conformal operator $\cB^0 \to \cB^1$ is isomorphic
to $ \mathbb{V}$; and the complex gives a resolution of this space
viewed as a sheaf.  The geometry of the manifold is partly encoded in
the coefficients of the differential operators in these sequences, and
in the conformally flat setting of the sphere the first PDE operator
$D_0:\cB^0\to \cB^1$ is fully integrable. From representation theory
(see \cite{BoeColl} and references therein) one can deduce that there
is a prolonged system and connection $\nabla^{D_0}$, equivalent to
$D_0$ in the sense discussed in Section \ref{prolongs}. In this setting
the space of solutions $\mathbb{V}$ is maximal in that $\dim
\mathbb{V}$ achieves the dimension bound for the maximal size of the
solution space. Since the sequence is a differential complex, the next
operator in the sequence $D_1$ gives an integrability condition for
$D_0$ and so on. That the complex gives a resolution means that these
integrability conditions are in a sense maximally severe.

It turns out that there are
``curved analogues'' of these sequences, these are the conformal cases
of the (generalised) Bernstein-Gelfand-Gelfand (BGG) sequences, a
class of sequences of differential operators that exist on any
parabolic geometry \cite{CSS,MikeNotices}. Unfortunately in the
general curved setting these sequences are no longer complexes, which
limits their applications.

As well as the operators $D_i:\cB^i\to \cB^{i+1}$ of the BGG sequence,
in even dimensions there are conformally invariant ``long operators''
$ \cB^k\to \cB^{n-k}$ for $ k=1,\cdots ,{n/2-1}$ \cite{BoeColl}.
Thus there are sequences of the form
\begin{equation}\label{gdet}
\cB^0\stackrel{D_0}{\to}\cB^1 \stackrel{D_1}{\to}\cdots 
\stackrel{D_{k-1}}{\to} \cB^k \stackrel{L_k}{\to} \cB^{n-k} 
\stackrel{D_{n-k}}{\to}\cdots \stackrel{D_{n-1}}{\to} \cB^n~.
\end{equation}
and, following \cite{BrGodeRham} (see also \cite{Tsrni}), we term
these detour sequences since, in comparison to the BGG sequence, the
long operator here bypasses the middle of the BGG sequence and the
operator $L_k$ takes us (or ``detours'') directly from $\cB^k$ to
$\cB^{n-k}$. Once again using the classification it follows
immediately that these detour sequences are in fact complexes in the
case that the structure is conformally flat. The operator $L_k$ is
again an integrability condition for $D_{k-1}$ but since $L_k$ has
higher order than the operator $D_k$ one does not expect that the
detour complex is a resolution. In fact, by for example considering
Taylor series expansions for solutions, it is straightforward to show
that there is local cohomology at the bundle $\cB^k$. We do not need
the details here, the main point is that $L_k$ is weaker, as an
integrability condition, than $D_k$.

\subsection{Curved detour complexes}\label{curveddets}
 An interesting direction is to 
try to find curved analogues of these complexes. That is detour
sequences that are actually complexes on conformally curved
structures. Remarkably this works in general in the de Rham case, that
is the case where the bundles $\cB^i$ are exterior powers of the
cotangent bundle, $\Lambda^i$. We write $\Lambda_k$ to denote the bundle of
conformal density weighted k-forms $\Lambda^k[2k-n]$; sections of this
pair conformally with $\Lambda^k$ in integrals. Then the formal
adjoint of the exterior derivative $d$ acts conformally $\delta:
\Lambda_i\to \Lambda_{i-1}$.
\begin{theorem} \cite{BrGodeRham}
In even dimensions there are conformally invariant differential operators 
$$
L_k: \Lambda^k \to \Lambda_k \quad \quad\quad L_k = (\delta d)^{n/2-k}+\mbox{ {\em lower order terms}}.
$$ 
so that 
$$
{\Lambda}^0\stackrel{d}{\to}\cdots\stackrel{d}{\to}{\Lambda}^{k-1}\stackrel{d}{\to}
{\Lambda}^k\stackrel{L_k}{\to}{\Lambda}_k\stackrel{\delta}{\to}{\Lambda}_{k-1}
\stackrel{\delta}{\to}
\cdots\stackrel{\delta}{\to}{\Lambda}_0
$$
is a conformal complex. In Riemannian signature the complex is elliptic.
\end{theorem}
\noindent The operator $L_{n/2-1}$ is the usual Maxwell operator
$\delta d$, while $L_0$, included here (cf. above) as an extreme case,
is the critical (i.e.\ dimension order) conformal Laplacian of GJMS
\cite{GJMS}. For $k\neq 0$ the complex may be viewed as a differential
form analogue of the critical GJMS operator. For
$0\leq k \leq n/2-2$ the $L_k$ have the form $\delta Q_{k+1} d$ where
the $Q_{\ell}$ are operators on closed forms that generalise the
$Q$-curvature; for example they give conformal pairings that descend
to pairings on de Rham cohomology \cite{newQ}.

\subsection{Variational constructions} \label{var}

On curved backgrounds, obtaining complexes for almost any other case
of \nn{gdet} seems, at first, to be hopeless. The point is this. One
of the main routes to constructing BGG sequences is via some variant
of the curved translation principle of Eastwood and others
\cite{ER,Baston,CSS,CD}. In the first step of this, one replaces the
de Rham complex with the sequence obtained by twisting this sequence
with an appropriate bundle and connection. In fact the connections
used are usually so called {\em normal conformal tractor connections}
as in \cite{CapGotrans}. These are connections on prolonged systems,
along the lines of the connection constructed explicitly in Section
\ref{prolongs}, but these are natural to the conformal structure and
are normalised in way that means that they are unique up to isomorphism.
Appropriate differential splitting operators are then used to extract
from the twisted de Rham sequence the sought BGG sequence. The full
details are not needed here, the main point for our current discussion is as follows.  If
the structure is conformally flat then the tractor connection is flat
and the twisted de Rham sequence sequence is still a complex. It
follows from this that the BGG sequence is a complex. However in
general the curvature of the tractor connection will obstruct the forming of a
complex. This is evident even at the initial stages of the sequence:
if we write $\Lambda^i(\mathcal{V})$ for the twistings of $i$-forms by
some tractor bundle then, writing $\nabla$ for the tractor connection,
the composition
\begin{equation}\label{comp}
{\Lambda}^0(\mathcal{V})\stackrel{d^\nabla}{\to}{\Lambda}^1 (\mathcal{V}) \stackrel{d^\nabla}{\to}{\Lambda}^{2}(\mathcal{V})
\end{equation}
is simply the curvature of $\nabla$ acting on $ \mathcal{V}$.

Ignoring this difficulty for the moment, an interesting case from
conformally flat structures is the conformal deformation complex of
\cite{GasqG}. Some notation first.  We write $\Lambda^{k,\ell}$ for
the space of trace-free covariant $(k+\ell)$-tensors $t_{a_1\cdots
a_kb_1\cdots b_\ell}$ which are skew on the indices $a_1\cdots a_k$
and also on the set $b_1\cdots b_\ell$. Skewing over more than $k$
indices annihilates $t$, as does symmetrising over any 3 indices.
Thus for example $\Lambda^{1,1}=S^2_0$, but in the case of this bundle
we will postpone changing to this notation.  In dimensions $n\geq 5$
the initial part of this BGG complex is,
\begin{equation}\label{flatdef}
\begin{picture}(350,20)(13,-4)
\put(70,0)
{$T \stackrel{{\sf K}_0}{\to}S^2_0[2]\stackrel{\sf C}{\to} 
\Lambda^{2,2}[2] \stackrel{\sf Bi}{\to}\Lambda^{3,2}[2] \to 
\cdots 
$}                                                         
\end{picture}
\end{equation} 
where $T$ is the tangent bundle.  Here ${\sf C}$ is the linearisation,
at a conformally flat structure, of the Weyl curvature as an operator
on conformal structure; ${\sf Bi}$ is a conformal integrability
condition arising from the Bianchi identity; the operator ${\sf K}_0$
is the conformal Killing operator, viz the operator which takes
infinitesimal diffeomorphisms to their action on conformal
structure. Since infinitesimal conformal variations take values in
$S^2_0[2]$, it follows from these interpretations that the cohomology
of the complex, at this bundle, may be interpreted as the formal
tangent space to the moduli space of conformally flat structures.

From this picture, and also from the discussion around \nn{comp}, we
do not expect a curved generalisation of the complex \nn{flatdef}.
One way to potentially avoid the issues brought up with \nn{comp} is
to consider detour sequences of the form \nn{gdet} with $k=1$; these
will be termed {\em short detour sequences}. (In \cite{GoPetob} there
is some discussion and applications of detour complexes for this BGG
in the conformally flat setting.) It turns out that this idea is
fruitful. One construction is based around the so-called obstruction
tensor $\mathcal{B}_{ab}$ of Fefferman and Graham \cite{FeffGrast}
which generalises to even dimensions the Bach tensor of Bach's 
gravity theory.  For our current purposes we need only to
know some key facts about this, and we shall introduce these as
required.  It is a trace-free conformal conformal 2-tensor with
leading term $\Delta^{n/2-2} C$ (where, recall, $C$ denotes the Weyl
curvature).  Let us write $K_0^*$ for the formal adjoint of $K_0$ and
$B$ for the linearisation, at the metric $g$, of the operator which
takes metrics $g$ to $\mathcal{B}^g$. By taking the Lie derivative of
$\cB_{ab}$, and using the fact \cite{GrHir} that $\cB_{ab}$ is the
total metric variation of an action (viz.\ $\int Q$ where $Q$ is
Branson's $Q$ curvature \cite{Brsharp}) we obtain the following (where
to simplify notation we have omitted the conformal weights).
\begin{theorem}\cite{BrGodefdet}\label{obflat} On even dimensional 
pseudo-Riemannian manifolds with the Fefferman-Graham obstruction
tensor vanishing everywhere, the sequence of operators
\begin{equation}\label{defdet}
 T\stackrel{K_0}{\to} S^2_0  \stackrel{B}{\to}
S^2_0 \stackrel{K_0^*}{\to} T
\end{equation}
is a formally self-adjoint complex of conformally invariant
operators. In Riemannian signature the complex is elliptic. 
\end{theorem}
\noindent As with the deformation complex \nn{flatdef} there is an
immediate interpretation of the cohomology at $S^2_0$. By construction
it is the formal tangent space to the moduli space of obstruction-flat
structures. There are determinant quantities, with interesting
conformal behaviour, for detour complexes \cite{BrGoprogress,Tsrni};
this detour torsion  should be especially interesting for
\nn{defdet}.

Leaving aside the potential applications it is already interesting that
there is a differential complex along these lines in such a general
setting. Although the obstruction tensor is rather mysterious, and at
this point has not been fully explored, it is known that it vanishes on
conformally Einstein manifolds \cite{FeffGrast,GoPetob,GrHir}, certain
products of Einstein manifolds \cite{GoLeit}, and on half-flat
structures in dimension 4.

\subsection{The Yang-Mills detour complex} \label{YM}
\newcommand{\nds}{{\tilde{\nabla}}}
\newcommand{\IT}[1]{{\rm(}{\it{\!#1}}{\rm)}} The construction in
Theorem \ref{obflat} and subsequent observations suggest the
possibility of a rich theory of (short) detour complexes. On a large
class of manifolds, the formally self-adjoint operator $B$ is
evidently sufficiently ``weak'' as an integrability condition for
$K_0$ that we obtain a complex. The construction may obviously be
generalised by using different Lagrangian densities (cf.\ $Q$) in the
first step.  On the other hand, if we use directly the
ideas from the {\em proof} of Theorem \ref{obflat}, then it seems that we 
may only obtain detour complexes
with either the Killing operator on vector fields,
$k\mapsto \mathcal{L}_k g$ (with $g$ the metric), or its conformal
analogue $K_0$, as the first operator $\cB^0\to \cB^1$.
 This motivates a rather different approach.

The simplest example of a conformal
detour complex is the Maxwell detour complex in dimension 4
\begin{equation}\label{Maxdet}
\Lambda^0 \stackrel{d}{\to} \Lambda^1 \stackrel{\delta d}{\to} 
\Lambda_1\stackrel{\delta}{\to}\Lambda_0~.
\end{equation} 
In other dimensions this is also a complex, but the middle operator is
no longer conformally invariant. Let us temporarily relax the
condition of conformal invariance and consider this complex on 
a pseudo-Riemannian $n$-manifold $(M,g)$ of signature
$(p,q)$, with  $n\geq 2$.

 For any vector vector bundle $V$, with
connection $\nds$, we might consider twisting the Maxwell detour. Of
course the result would not be a complex. However by a minor variation
on this theme we make an interesting observation.
Recall that we write $d^\nds$ for
the connection-coupled exterior derivative operator
$d^\nds:\Lambda^k(V)\to \Lambda^{k+1}(V)$. Of course we could equally
consider the coupled exterior derivative operator
$d^{\nds}:\Lambda^k(V^*)\to \Lambda^{k+1}(V^*)$ and for the
formal adjoint of this we write
$\delta^\nds:\Lambda_{k+1}(V)\to\Lambda_k(V)$.

Denote by $F$ the curvature of $\nds$ and write $F\cdot$ for the action
of the curvature on the twisted 1-forms, $F\cdot :\Lambda^1(V)\to
\Lambda_1(V)$ given by
$$ 
(F\cdot \phi)_a:= F_a{}^b\phi_b, 
$$ where we have indicated the abstract form indices explicitly, whereas
the standard $\End(V)$ action of the curvature on the $V$-valued
1-form is implicit. Using this we construct a
differential operator
$$
M^{\nds}:\Lambda^1(V)\to \Lambda_1(V)
$$
by 
$$
M^{\nds}\phi=\delta^\nds d^\nds\phi-F\cdot \phi .
$$ By a direct calculation, the composition $M^\nds d^\nds $ amounts to an exterior
algebraic action by the ``Yang-Mills current'' $\delta^\nds F$, thus
we have the first statement of the following result. The other claims
are also easily verified.
\begin{theorem}\cite{GSS}\label{twistthm}
The sequence of operators,
\begin{equation}\label{detseq}
\Lambda^0(V)\stackrel{d^\nds}{\to} \Lambda^1(V)\stackrel{M^\nds}{\to} 
\Lambda_1(V)\stackrel{\delta^\nds}{\to}\Lambda_0(V)
\end{equation}
is a complex if and only if 
the curvature $F$ of the connection $\nds$ satisfies the (pure) Yang-Mills
equation
$$
\delta^\nds F=0.
$$ In addition:\\ \IT{i} If $\nds$ is an orthogonal or unitary
connection then the sequence is formally self-adjoint.\\ \IT{ii} In
Riemannian signature the sequence is elliptic.\\
\IT{iii} In dimension 4 the complex is conformally invariant.
\end{theorem}

This obviously yields a huge class of complexes. For example, taking
$V$ to be any tensor (or spin) bundle on a (spin) manifold with
harmonic curvature (i.e. a pseudo-Riemannian structure where the
Riemannian curvature satisfies the Yang-Mills equations) yields a
complex.  Einstein metrics are harmonic so this is large class of
structures.  To be specific if we take, for example, $V$ to be the second exterior power of the tangent bundle
$T^2$, and use $\nabla$ to denote the Levi-Civita connection
then we obtain the complex 
$$ T^2 \stackrel{d^\nabla}{\to} \Lambda^1\otimes T^2
\stackrel{M^\nabla}{\to} \Lambda^1\otimes T^2
\stackrel{\delta^\nabla}{\to} T^2 , 
$$
where $M^\nabla$ is given by
$$
S_b{}^{cd}\mapsto 
-2\nabla^a\nabla_{[a } S_{b]}{}^{cd}- R_{ba}{}^c{}_eS^{aed} - R_{ba}{}^d{}_eS^{ace}~.
$$ This example reminds us that part (iii) of the theorem means that
the complex \nn{detseq} is conformally invariant in dimension 4,
provided the connection $\nds$ is conformally invariant. For
applications in conformal geometry we need suitable conformal
$(V,\tilde{\nabla})$.

\subsection{Short detour complexes and tractor connections}\label{diagram}

\newcommand{\cD}{D}
\renewcommand{\d}{\delta}

 We may use Theorem \ref{twistthm} to construct more 
differential complexes. Consider the following general situation. Suppose 
that there are vector
bundles (or rather section spaces thereof) $\cB^0$, $\cB^1$, $\cB_1$ and
$\cB_0$ and differential operators $L_0$, $L_1$, $L^1$, $L^0$, $\cD$
and $\overline{\cD}$ which act as indicated in the following diagram:
 \newcommand{\fbbbb} {\mbox{$
\begin{picture}(40,2)(1.6,-2)
\put(1,1.0){\line(1,0){58}}
\end{picture}$}} 
$$
\begin{picture}(350,100)(12,-80)

\put(10,-25){[D]}

\put(70,0) {$ \Lambda^0(V)\hspace*{2pt}\stackrel{d^\nds}{\longrightarrow}
\hspace*{8pt} \Lambda^1(V) \stackrel{M^\nds}{\fbbbb\longrightarrow}
\hspace*{8pt} 
\Lambda_1(V)\hspace*{2pt} \stackrel{\d^\nds}{\longrightarrow}\hspace*{8pt} 
\Lambda_0(V)
$}                                                         

\put(77,-38){\vector(0,1){31}}                               
\put(66,-25){\scriptsize$L_0$}

\put(143,-38){\vector(0,1){31}}                               
\put(132,-25){\scriptsize$L_1$}

\put(249,-6){\vector(0,-1){31}}                              
\put(238,-25){\scriptsize$L^1$}

\put(315,-6){\vector(0,-1){31}}                              
\put(304,-25){\scriptsize$L^0$}

\put(70,-50)
{$
\hspace*{2pt} \cB^0{\phantom{(V)}} \stackrel{\cD}{\longrightarrow} 
\hspace*{8pt} \cB^1{\phantom{(V)}} \hspace*{-2pt}
\stackrel{M^{\cB}}{\fbbbb\longrightarrow} \hspace*{8pt}
\cB_1{\phantom{(V)}}
\stackrel{\overline{\cD}}{\longrightarrow}\hspace*{8pt}  \cB_0{\phantom{(V)}}
$}                                                         

\end{picture}
$$
The top sequence is \nn{detseq} for a connection $\nds$ with curvature
$F$ and the operator $M^\cB:\cB^1\to \cB_1$ is defined to be the
composition $L^1M^{\nds}L_1$.  Suppose that the squares at each end
commute, in the sense that as operators $\cB^0\to \Lambda^1(V)$ we have
$d^{\nds}L_0=L_1\cD$ and as operators $\Lambda_1(V)\to \cB_0$ we have
$L^0\d^\nds=\overline{\cD}L^1$.  Then on $\cB^0$ we have
$$
M^{\cB} \cD = L^1 M^{\tilde{\nabla}} L_1 \cD = L^1 M^{\tilde{\nabla}} d^{\tilde{\nabla}} L_0 = L^1\epsilon(\d^{\tilde{\nabla}} F) L_0 ,
$$ and similarly $\overline{\cD} M^{\cB}=- L^0 \i(\d^{\tilde{\nabla}}
F) L_0 $. Here $\epsilon(\d^{\tilde{\nabla}} F)$ and
$\i(\d^{\tilde{\nabla}} F)$ indicate exterior and interior actions of
the Yang-Mills current $\d^{\tilde{\nabla}} F$.  Thus if
$\tilde{\nabla}$ is Yang-Mills then the lower sequence, viz.
\begin{equation}\label{transcx}
\cB^0 \stackrel{\cD}{\longrightarrow}  \cB^1
\stackrel{M^{\cB}}{\longrightarrow} 
\cB_1
\stackrel{\overline{\cD}}{\longrightarrow}  \cB_0 ~,
\end{equation}
is a complex.

\noindent {\bf Remark:} Note that if the connection $\nds$ preserves a
Hermitian or metric structure on $V$ then we need only the single commuting 
square $d^{\nds}L_0=L_1\cD$ on $\cB^0$ to obtain such a complex; by taking formal adjoints we obtain a 
second commuting square $(L^0\d^\nds=\overline{\cD}L^1): \cB_1\to \cB_0$ where 
$\cB_0$ and $\cB_1$ are appropriate density twistings of the bundles 
$\cB^0$ and $\cB^1$ respectively.
\quad \endrk 

We are now in a position to link the various constructions and
indicate the general prospects. Suppose that in the setting of
Riemannian or conformal geometry, we have some finite type PDE
operator $B:\cB^0\to \cB^1$ that we wish to study. A question we
raised in Section \ref{prolongs} is whether there is a corresponding
connection $\nabla^D$ that is equivalent in the sense that solutions
of $D$ are in 1-1 correspondence with parallel sections of $\nabla^D$.
Consider the situation of the first square in the diagram. If $ D \phi
=0$ then obviously $L_0 \phi$ is parallel for $\tilde{\nabla}$. Now
suppose that the $L_i$ are in fact differential splitting operators,
that is for $i=0,1$ these may be inverted by differential operators
$L_i^{-1}$ in the sense that $L_i^{-1}\circ L_i=\id_{\cB^i}$. Then $D=
L_1^{-1}(L_1 D)= L_1^{-1}(\tilde{\nabla}L_0)$ and so solutions for $D$
are mapped by $L_0$ injectively to sections of $V$ that are parallel
for $\tilde{\nabla}$. Finally if also $L^{-1}_0$ is well defined on
general sections of $V$ and $L_0\circ L_0^{-1}$ acts as the identity
on the space of parallel sections of $V$, then $L^{-1}_0$ maps
parallel sections for $\tilde{\nabla}$ to $D$-solutions. All these
conditions are somewhat more than is strictly necessary according to
ideas of section \ref{prolongs}. Nevertheless, as we indicate below,
these are satisfied for some examples, and it seems likely that these
are the simplest cases of a large class.

\subsection{Examples}

\newcommand{\J}{J}
\newcommand{\Rho}{P}

Here we treat mainly a single example (with a hint of a second example), but
it illustrates well the general idea. We return to the setting of
conformal $n$-manifolds.  Modulo the trace part, the conformal
transformation of the Schouten tensor is controlled by the equation
\begin{equation}\label{sc}
D\si=0 \quad \mbox{where}\quad D\si :={\rm TF}(\nabla_{a} \nabla_{b}\si+P_{ab}\si ),
\end{equation}
which is written in terms of some metric $g$ from the conformal class.
In particular a metric $\si^{-2}\bg$ is Einstein if and only if the
scale $\si\in \ce[1]$ is non-vanishing and satisfies \nn{sc}.

The conformal standard tractor bundle and connection arises from a
prolongation of this overdetermined equation as follows, cf.\
\cite{BEGo,CapGoluminy}.  Recalling the notation developed in Section
\ref{var} we have $\Lambda^{1,1}$ as an alternative notation for
$S^2_0$.  Let us also write $\Lambda_{1,1}:=\Lambda^{1,1}\otimes
\ce[4-n]$. As usual we use the same notation for the spaces of
sections of these.  The {\em standard tractor bundle} $\bT$ may be
defined as the quotient of $J^2 \ce[1]$ by the image of
$\Lambda^{1,1}[1]$ in $J^2 \ce[1]$ through the jet exact sequence at
2-jets. Note that there is a tautological operator $\bD:\ce[1]\to
\Lambda^0(\bT)$ which is simply the composition of the universal 2-jet
operator $j^2:\ce[1]\to J^2\ce[1]$ followed by the canonical
projection $J^2\ce[1]\to \Lambda^0(\bT)$. By construction both $\bT$
and $\bD$ are invariant, they depend only on the conformal structure
and no other choices.  Via a choice of metric $g$, and the Levi-Civita
connection it determines, we obtain a differential operator 
\begin{equation}\label{Dsplit}
  \ce[1]\to \ce[1]\oplus \Lambda^1[1]\oplus \ce[-1] \quad \mbox{by} \quad \si\mapsto ( \si,
  \nabla_a \si, - \frac{1}{n}(\Delta + \J) \si ).
\end{equation} 
Since this is second order it factors through a linear map
$J^2\ce[1]\to\ce[1]\oplus \Lambda^1[1]\oplus \ce[-1] $. Considering
Taylor series one sees that the kernel of this is a copy of $\Lambda^{1,1}[1]$
and so the map determines
an isomorphism
\begin{equation}\label{split}
\bT \stackrel{g}{\cong} \ce[1]\oplus \Lambda^1[1]\oplus
\ce[-1] ~.
\end{equation}
In terms of this the formula at the right extreme of the display \nn{Dsplit} 
then tautologically gives an
explicit formula for $\bD$.  
This is a differential splitting operator since through
the jet projections there is conformally invariant surjection $X:
\ce(\bT)\to \ce[1]$ which inverts $\bD$. In terms of the splitting
\nn{split} this is simply $(\si,\mu,\rho)\mapsto \si$.

Observe that if we change to a conformally related metric
$\widehat{g}=e^{2\omega}g$ then the  Levi-Civita connection has a
corresponding transformation, and so we obtain a different isomorphism.
For $t\in \bT$ then via \nn{split} we have $[t]_g=(\si,\mu,\rho)$. In
terms of the analogous isomorphism for $\widehat{g}$ we have
\begin{equation}\label{trans}
[t]_{\widehat{g}} = (\widehat{\si},\widehat{\mu}_a,\widehat{\rho})=(\si,\mu_a+\si\Upsilon_a,\rho-\bg^{bc}\Upsilon_b\mu_c-
\tfrac{1}{2}\si\bg^{bc}\Upsilon_b\Upsilon_c). 
\end{equation}
where $\Upsilon=d\omega$.

Now let us
 define a connection on
$\ce[1]\oplus \Lambda^1[1]\oplus \ce[-1] $ by the formula
\begin{equation}\label{trconn}
\nabla_a
\left( \begin{array}{c}
\si\\\mu_b\\ \rho
\end{array} \right) : =
\left( \begin{array}{c}
 \nabla_a \si-\mu_a \\
 \nabla_a \mu_b+ g_{ab} \rho +\Rho_{ab}\si \\
 \nabla_a \rho - \Rho_{ab}\mu^b  \end{array} \right) 
\end{equation}
where, on the right-hand-side $\nabla$ is the Levi-Civita connection
for $g$.  Obviously this determines a connection on $\bT$ via the
isomorphism \nn{split}. What is more surprising is that if we repeat
this using a different metric $\widehat{g}$, this induces the {\em
  same} connection on $\bT$. Equivalently the connection in the
display transforms according to \nn{trans}. This canonical connection
on $\bT$ depends only on the conformal structure and is known as the
{\em (standard) tractor connection}.

There are several ways we may understand the connection \nn{trconn}
and its invariance. On the one hand this may be viewed as a special
case of a {\em normal} tractor connection for parabolic geometries. This construction is treated from that point of view in
\cite{CapGoluminy}.  Such connections determine the normal Cartan
connections on the corresponding adapted frame bundles for the given
tractor bundle $\mathbb{T}$, \cite{CapGotrans}.  On the other hand one
may start directly with the operator $D$ in \nn{sc}. It is readily
verified that, for $n\geq 3$, this is of finite type and so we may
attempt to construct a prolonged system and connection $\nabla^D$ 
for the operator along the
lines of the treatment of Killing operator in Section \ref{prolongs}.
This is done in detail in \cite{BEGo} and one obtains exactly \nn{trconn}.  
So the connection \nn{trconn} is the sought $\nabla^D$.
Some key points are as follows: 
since the operator in \nn{sc} acts $D:\ce[1]\to \Lambda^{1,1}[1]$ and is second
order,  it  factors through a linear bundle map
$D^{(0)}:J^2\ce[1]\to\Lambda^{1,1}[1] $; from the formula for $D$ it
is clear that this is a splitting of the exact sequence
$$
0\to \Lambda^{1,1}[1]\to J^2\ce[1]\to \bT\to 0
$$
that defines $\bT$. It follows that the prolonged system will include
$\mathcal{T}$. In this case the key to the closing up of the prolonged
system is that any completely symmetric covariant 3 tensor that is
also pure trace on some pair of indices is necessarily zero. From this it 
follows that the full prolonged system may be expressed in terms of  $\bT$.

There is also a
differential splitting operator 
$$
 E:\Lambda^{1,1}[1] \to \Lambda^1(\mathcal{T}) \quad \quad \quad \psi_{ab} \mapsto
(0,\psi_{ab},-(n-1)^{-1}\nabla^b\psi_{ab})
$$ (cf.\ \cite{Esrniconf}). An easy calculation verifies this is
also conformally invariant, and from the formula it is manifest that
we have  an operator $E^{-1}$ (in fact bundle map) that inverts
$E$ from the left. Crucially, we have the following.
\begin{proposition} \cite{GSS} \label{eincomm}
 As differential operators on $\ce[1]$, we have
$$
\nabla^D \bD = E D   ~.
$$ For $\si\in \ce[1]$, $\bD\si$
is parallel if and only if $D\si=0$.
\end{proposition}
\noindent{\bf Proof:} The second statement is immediate from the
first.  A straightforward calculation verifies that either composition
applied to $\si\in \ce[1]$ yields
$$
\left( \begin{array}{c}
0\\ {\rm TF}(\nabla_{a} \nabla_{b}\si+P_{ab}\si ) \\ 
-\frac{1}{n}\nabla_a(\Delta \si +J \si)-P_a{}^c\nabla_c\si
\end{array} \right) 
$$ \quad $\blacksquare$ \\ 

For our present purposes the main point of Proposition \ref{eincomm}
is that it gives the first step in constructing a special case of a
commutative detour diagram [D]. First note that we have exactly the
situation discussed in the last part of Section \ref{diagram}: here
$\bD$ and $E$ play the roles of, respectively, $L_0$ and $L_1$. They
are differential splitting operators and have inverses exactly as
discussed there. It is an easy exercise to verify that if $I$ is a
parallel section of $\bT$ then $I=D\si$ for some section $\si$ of
$\ce[1]$, as observed in \cite{GoNur} (and so $X$ as an inverse to
$\bD$ maps parallel tractors to solutions of $D$). A useful consequence
of this last result is that a conformal manifold with a parallel tractor is
{\em almost Einstein} in the sense that it has a section of $\ce[1]$
that gives an an Einstein scale on an open dense subset (see
\cite{goalmost} for further details).
  
Next we observe that  there is a
conformally invariant {\em tractor metric} $h$ on $\mathcal{T}$ given
(as a quadratic form) by $(\si,~\mu,~\rho)\mapsto
\bg^{-1}(\mu,\mu)+2\si \rho$. This 
has signature $(p+1,q+1)$ (corresponding to $\bg$ of signature
$(p,q)$) and is preserved by the tractor connection.  
In view of this, and using our general observations in Section 
\ref{diagram},  the formal
adjoints of the operators above give the other required commutative square of
operators. That is with
$$\begin{array}{ccll}
 && D^*: \Lambda_{1,1}[-1]\to \Lambda_0[-1] & 
\phi_{ab}\mapsto \nabla^{a}\nabla^b \phi_{ab}+ P^{ab}\phi_{ab}  \\ 
&& E^*: \Lambda_1(\mathcal{T}) \to \Lambda_{1,1}[-1] & (\alpha_a,~\nu_{ab},~\tau_a) \mapsto \nu_{(ab)_0}
+\frac{1}{n-1}\nabla_{(a}\alpha_{b)_0} \\
 && \bD^*: \Lambda_0(\mathcal{T}) \to \Lambda_0[-1]& (\si,~\mu_b,~\rho) \mapsto \rho -\nabla^a\mu_a-\frac{1}{n}(\Delta \si+ J\si)\\
&&  \d^{\nabla^D}: \Lambda_{1}(\mathcal{T})\to \Lambda_0(\mathcal{T}) & \Phi_{aB}\mapsto -\nabla^a \Phi_{aB}
\end{array}
$$  
we have $\bD^*\delta^{\nabla^D}= D^* E^*$ on $\Lambda_{1}(\mathcal{T})$.

We now want to consider what it means for the tractor connection to
satisfy the Yang-Mills equations.  
 Using the explicit presentation
\nn{trconn} of the tractor connection at a metric, it is straightforward to
calculate its curvature. This is
$$
\Om_{ab}{}^C{}_D=
\left(\begin{array}{lll}
0&0&0\\
A^{c}{}_{ab} & C_{ab}{}^c{}_d & 0\\
0 & -A_{dab} & 0
\end{array}\right)
$$ where $A$ is the {\em Cotton tensor}, $ A_{abc}:=2\nabla_{[b}P_{c]a}$.
Taking the required divergence we obtain $-\delta^{\nabla^{D}}\Omega$
(see e.g.\ \cite{GoNur} for further details),  
\begin{equation}\label{divtr}
\nabla^D_a \Om^a{}_{b}{}^C{}_E= \left(\begin{array}{lll}
0&0&0\\
B^{c}{}_{b} & (n-4)A_{b}{}^c{}_e & 0\\
0 & -B_{eb} & 0
\end{array}\right)
\end{equation}
 where, on the left-hand side, $\nabla^D$ is really the Levi-Civita connection
coupled with the tractor connection on $\End(\mathcal{T})$ induced from
$\nabla^D$.  Here $ B_{ab}$ is the 
{\em Bach tensor} $ B_{ab}:=\nabla^c
A_{acb}+P^{dc}C_{dacb}$. 

Let us say that a {\em pseudo-Riemannian} manifold is
{\em semi-harmonic} 
 if its tractor curvature is Yang-Mills, that
is $ \nabla^a \Om_{ab}{}^C{}_D=0$. Note that in dimensions $n\neq 4$ this
is not a conformally invariant condition and a semi-harmonic
space is a Cotton space that is also Bach-flat.   From our 
observations above, the
semi-harmonic condition is conformally invariant in dimension 4
and according to the last display we have the following
result (which in one form or another has been known for many years e.g.\ 
\cite{Merk,Lew}). 
\begin{proposition}\label{trYM}
 In dimension 4 the conformal tractor connection is a Yang-Mills
connection if and only if the structure is Bach-flat.
\end{proposition}
 Writing $M^{\mathcal{T}}$ for the composition $E^* M^\nabla E $
 from the construction in section \ref{diagram} above, we have the
 following (except for the claim of ellipticity, which is straightforward).
\begin{theorem}\cite{GSS}\label{paracase}
The sequence
\begin{equation}\label{Einseq}
 \Lambda^0[1]\stackrel{P}{\to} \Lambda^{1,1}[1] 
\stackrel{M^{\mathcal{T}}}{\longrightarrow} \Lambda_{1,1}[-1]
\stackrel{P^*}{\to} \Lambda_0[-1]
\end{equation}
has the following properties. \\
\IT{i} It is a formally self-adjoint sequence of
differential operators and, 
for $\si\in \ce[1]$ 
\begin{equation}\label{MP}
(M^{\mathcal{T}} D \si)_{ab} =
-TFS\big( B_{ab}\si - (n-4) A_{abc}\nabla^c\si \big) ,
\end{equation}
where $TFS(\cdots)$ indicates the trace-free symmetric part of the
tensor concerned.  In particular it is a complex on
semi-harmonic manifolds. \\ 
\IT{ii} In the case of Riemannian
signature the complex is elliptic. \\ 
\IT{iii} In dimension 4,
\nn{Einseq} is sequence of conformally invariant operators and it is a
complex if and only if the conformal structure is Bach-flat.
\end{theorem}

\begin{corollary}\label{EimpliesB}
  Einstein 4-manifolds are Bach-flat. 
 \end{corollary}
\noindent{\bf Proof:} If a non-vanishing density $\si$ is an
Einstein scale then, calculating in that scale, we have $M^{\mathcal{T}} D\si
= -B \si$, where $B$ is the Bach tensor. On the other hand if $\si$ is
an Einstein scale then $D \si =0$ (see \nn{sc}).  \quad $\blacksquare$\\
In fact, the result in the Corollary, and more general results, have been 
known by other means for some
time (see e.g.\ \cite{GoPetob,GrHir} and references
therein). Nevertheless it seems the detour complex gives an
interesting route to this.

\vspace{4mm}
\renewcommand{\S}{{\Bbb S}}
We conclude here with just the statements of another example from
\cite{GSS}.  This is also obtained from the diagram [D] and an
appropriate tractor connection $\nabla$ for an operator $D$, in this
case the operator $D$ is the twistor operator. Thus 
we assume here that we have a conformal spin structure.
Following \cite{TomLeb}  we
 write $\S$ for the basic spinor bundle and $\overline{\S}={\S}[-n]$
 (i.e.\ the bundle that pairs globally in an invariant way with $\S$
 on conformal $n$-manifolds). The weight conventions here
 give $\S$ a ``neutral weight''. In terms of, for example, the Penrose
 weight conventions
 $\S=E^\lambda[-\frac{1}{2}]=E_{\lambda}[\frac{1}{2}]$, where
 $E^\lambda$ denotes the basic contravariant spinor bundle in
 \cite{otnt}.

 We write ${\rm Tw}$
for the so-called twistor bundle, that is the subbundle of
$\Lambda^1 \otimes \S[1/2]$ consisting of form spinors $u_a$ such that
$\gamma^a u_a=0$, where $\gamma_a$ is the usual Clifford symbol.  We
use $\S$ and ${\rm Tw}$ also for the section spaces of these bundles.
The {\em twistor operator} is the conformally invariant Stein-Weiss
gradient
$$
{\bf T}: \S[1/2]\to {\rm Tw}
$$
given explicitly by
$$
\psi\mapsto \nabla_a \psi +\frac{1}{n}\gamma_a\gamma^b \nabla_b \psi ~.
$$ 
This completes to a
differential complex as follows.
\begin{theorem}\cite{GSS}\label{twcase}
On semi-harmonic pseudo-Riemannian $n$-manifolds $n\geq 4$ we have a
differential complex
\begin{equation}\label{spcx}
 \S[1/2]  \stackrel{\bf T}{\to} {\rm Tw} 
\stackrel{\bf N}{\longrightarrow} \overline{\rm Tw}
\stackrel{{\bf T}^*}{\to} \overline{\S}[-1/2] ,
\end{equation} 
 where ${\bf T}$ is the usual twistor operator, ${\bf T}^*$ its
formal adjoint, and ${\bf N}$ is third order. The sequence is formally self-adjoint and in the case
of Riemannian signature the complex is elliptic.

In dimension 4 the sequence \nn{spcx} is conformally invariant and it
is a complex if and only if the conformal structure is Bach-flat.
\end{theorem}
\noindent The operator ${\bf N}$ is a  third order analogue of a  
Rarita-Schwinger operator.
Of course on a fixed pseudo-Riemannian
manifold we may ignore the conformal weights. 

\subsection{Outlook} 
BGG complexes do not generally have curved analogues. The weaker
integrability conditions involved with detour complexes suggest the
possibility of large classes of such objects. 

In the terminology of physics, the formally self-adjoint short detour
complexes $\cB^0\to \cB^1\to \cB_1\to \cB_0$ give ``classically
consistent systems''. The constraint equations, $\cB_1\to \cB_0$ which
give an integrability condition on the middle operator, are suitably
dual to, the gauge transformations of $\cB^1$ given by $\cB^0\to
\cB^1$. Thus if we take $\cB^1\to\cB_1$ as the field equations then
the first cohomology of the complex gives ``observable'' field
quantities. It seems that there is some scope for quantising this
picture \cite{GoW}.

As we see in Corollary \ref{EimpliesB} the first operator $\cB^0\to
\cB^1$ can  itself have an important geometric or physical
interpretation and the complexes provide a tool for studying these.
The first cohomology of the complex gives a global invariant.

Generalising the examples given above requires two main steps. The
first is to develop the theory of prolonged differential systems in a way
which leads to prolonged systems that share the invariance properties
of the original equation (this is currently the subject of joint work
with M.G. Eastwood).  The second part is to understand what operators
may be used to replace $M^{\nabla^D}$, so that, for example, we may
have conformally invariant examples in higher dimensions.

\end{document}